\documentstyle[12pt]{article}
\input epsf
\newtheorem{theorem}{Theorem}
\newtheorem{lemma}{Lemma}
\newtheorem{proposition}{Proposition}
\newtheorem{remark}{Remark}
\newtheorem{corollary}{Corollary}
\newtheorem{definition}{Definition}

\newlength{\rig}

\newlength{\hei}

\newcommand{\fgr}[3]{
\setlength{\rig}{0.21\textwidth} \setlength{\hei}{7.8cm}
\begin{figure}%[p]
\rule{\rig}{0in} \epsfxsize=9cm \epsffile{#1} \caption{#3}\label{#2}
\end{figure}
}

\newcommand{\text}[1]{{\mbox {#1}}}

\begin{document}

\author{{\bf Anton Savin\footnotemark \rule{5pt}{0pt} and \addtocounter{footnote}{-1}
Boris Sternin\footnotemark }\addtocounter{footnote}{-1}\thanks{
Supported by the Russian Foundation for Basic Research under grants N
99-01-01254, 99-01-01100 and 00-01-00161, 01-01-06013 and by
Institut f\"ur Mathematik, Universit\"at Potsdam.} \\
%EndAName
Moscow State University\\
[3mm] e-mail: antonsavin@mtu-net.ru\\
[3mm] e-mail: sternine@mtu-net.ru}
\title{{\bf Eta-invariant and Pontryagin duality in $K$-theory}}
\date{}
\maketitle

\begin{abstract}
The topological significance of the spectral Atiyah--Patodi--Singer $\eta$%
-invariant is investigated. We show that twice the fractional part of the
invariant is computed by the linking pairing in $K$-theory with the
orientation bundle of the manifold. The Pontryagin duality implies the
nondegeneracy of the linking form. An example of a nontrivial fractional
part for an even-order operator is presented.
\end{abstract}

\setcounter{page}{1}

\vspace{15mm}

\vspace{15mm}

{$\!\!\!$\bf Keywords}:$\!$ eta-invariant, $\!K$-theory, Pontryagin duality, linking
coefficients, Atiyah--Patodi--Singer theory, modulo n index

\vspace{15mm}

{\bf 2000 AMS classification}: Primary 58J28, 19L64, Secondary 58J22, 19K56,
58J40

\vfill

\newpage

\tableofcontents
\vspace{0.3cm}

\section*{Introduction}

\addcontentsline{toc}{section}{Introduction}

P. Gilkey noticed \cite{Gil7} that the Atiyah--Patodi--Singer $\eta $%
-invariant \cite{APS0} is rigid in the class of differential operators on
closed manifolds when the following condition is satisfied
\[
{\rm {ord}}A+\dim M\equiv 1\left( {\rm mod} \,2\right) .
\]
More precisely, in this case the fractional part of the spectral $\eta$%
-invariant of an elliptic self-adjoint differential operator $A$ is
determined by the principal symbol and it is a homotopy invariant of the
principal symbol.

In this framework the $\eta $-invariant of geometric first-order operators
on even-dimen\-sio\-nal manifolds was studied in \cite{Gil8}. It turned out that
the invariant takes only dyadic values. Moreover, nontrivial
fractional parts not equal to $1/2$ appear only on nonorientable manifolds.
As a typical example we should point out (self-adjoint) Dirac operator on
manifolds with a $pin^c$ structure in the tangent bundle. The operators of
this form are defined, for example, on the projective spaces ${RP}^{2n}.$ In
this case, the fractional part of the $\eta $-invariant is an important $%
pin^c$-cobordism invariant \cite{BaGi1}. The case of general differential
operators remained unexplored. Even the problem of nontriviality of the $\eta
$-invariant's fractional part for even-order operators on odd-dimensional
manifolds remained open. We answer these questions in the present paper.

The following natural question appears, when we try to compute a fractional
homotopy invariant: what terms could be used to express such an invariant?
One of the candidates is the expression of the invariant in terms of some
{\em linking indices\/} \cite{Pon2,Pon3}. These indices define the bilinear
pairing
\[
{\rm {Tor}}H_i\left( M\right) \times {\rm {Tor}}H_{n-i-1}\left( M\right)
\longrightarrow {Q}/{Z}
\]
of torsion homology classes on a closed oriented $n$-dimensional manifold.
The linking indices have numerous applications (see \cite{SeTr1},\cite{FoFu1}%
,\cite{Mil1},\cite{BrMo1},\cite{KeMi1}).

In the present paper, we study a similar linking form in $K$-theory and
elliptic theory (the latter appeared also in the paper \cite{MoWi1}). The
Pontryagin duality in $K$-theory permits us to prove the nondegeneracy of
the linking form. The main result of the paper is the equality of fractional
part of twice the $\eta$-invariant with the linking index of the principal
symbol of the elliptic operator with the orientation bundle of the manifold.
Using this equality and the properties of the linking form, we prove the
nontriviality of the fractional part of the $\eta $-invariant for some even
order operators.

The proofs of the above results are based on a reduction of the spectral
invariant under consideration to a homotopy invariant dimension functional
of subspaces defined by pseudodifferential projections. This functional was
introduced in \cite{SaSt1,SaSt2}. The index formula for elliptic operators
which act in such subspaces makes it possible to express the fractional part
of the dimension functional in terms of the index of an operator in elliptic
theory with coefficients in a finite cyclic group ${Z}_n$. Together with the
corresponding theorem for the index modulo $n$ this expresses the fractional
invariant under consideration in topological terms.

Let us briefly describe the contents of the paper. In the first section,
we recall  for
the reader`s convenience the relationship of the $\eta $-invariant
to the dimension functional of subspaces defined by pseudodifferential
projections. Here we briefly present the necessary results from \cite
{SaSt1,SaSt2,SaScS5}. We also express the fractional parts of the invariants
in terms of the index of some elliptic operators in subspaces. At the
beginning of the second section we describe the main results of \cite{SaScS3}
concerning index theory modulo $n$. The rest of the section is devoted to
the proof of the Pontryagin duality in $K$-theory. The
linking pairing is introduced in Section 3. The fourth section contains a computation of
the action of antipodal involution on $K$-groups of real vector bundles.
Here we use the description \cite{Karo1} of these groups in terms of
Clifford algebras. In the following section we state the main theorem,
which expresses the fractional part of twice the dimension functional in terms of
a linking index. Examples are presented. In particular, we construct an even
subspace on which the dimension functional takes a value with a nontrivial
fractional part by means of the ``cross product'' \cite{AtSi1} of elliptic
operators. This gives a positive answer to the question of P.~Gilkey
concerning the existence of even-order operators with nontrivial
fractional part of the $\eta $-invariant.

We would like to thank Prof. P.~Gilkey for a fruitful discussion
concerning the results obtained in this paper.

\noindent

\section{Eta-invariant and index in subspaces}

{\bf 1.} Let $E$ be a vector bundle on a smooth manifold $M.$ A linear
subspace $\widehat{L}\subset C^\infty \left( M,E\right) $ is called {\em %
pseudodifferential}, if it is the range of a pseudodifferential projection $P
$ of order zero
\[
\widehat{L}={\rm Im} \,P,\quad P:C^\infty \left( M,E\right) \rightarrow
C^\infty \left( M,E\right) .
\]
The vector bundle $L={\rm Im} \,\sigma \left(P\right) \subset \pi ^{*}E$
over the cospheres $S^{*}M$ is called the {\em symbol} of a subspace.

Consider the antipodal involution
\[
\alpha :T^{*}M\longrightarrow T^{*}M,\quad \alpha \left( x,\xi \right)
=\left( x,-\xi \right)
\]
of the cotangent bundle $T^*M$. A subspace $\widehat{L}\subset C^\infty
\left( M,E\right) $ is {\em even (odd)}, if its symbol $L$ is invariant
(antiinvariant) under the involution:
\begin{equation}
L=\alpha ^{*}L,\quad {\mbox {{or }}}{\quad }L\oplus \alpha ^{*}L=\pi ^{*}E.
\label{omg}
\end{equation}
Denote by $\widehat{{\rm {Even}}}\left( M\right) $ $\left( \widehat{{\rm {Odd%
}}}\left( M\right) \right) $ the semigroups of all even (odd)
pseudodifferential subspaces with respect to the direct sum.

Pseudodifferential subspaces can be also defined by means of self-adjoint
elliptic operators. In this case the parity of the subspace corresponds to
the parity of order for differential operators.

\begin{proposition}
Let $A$ be an elliptic self-adjoint operator of a nonnegative order. Then
the subspace $\widehat{L}_{+}(A)$ generated by eigenvectors of $A$, which
correspond to nonnegative eigenvalues, is pseudodifferential and its symbol $%
L_{+}\left( A\right) $ is equal to the nonnegative spectral subbundle of the
principal symbol $\sigma \left( A\right) :$%
\[
L_{+}\left( A\right) =L_{+}\left( \sigma \left( A\right) \right) \in {\rm %
Vect}(S^{*}M).
\]
Suppose that $A$ is a differential operator or, more generally, it satisfies
the following condition
\[
\alpha ^{*}\sigma \left( A\right) =\pm \sigma \left( A\right)
\]
on its principal symbol. Then the subspace $\widehat{L}_{+}(A)$ is even
(odd).
\end{proposition}

The proof of this proposition can be found in \cite{BiSo1} or \cite{SaSt2}.

Even subspaces on odd-dimensional manifolds and odd subspaces on
even-dimensio\-nal ones admit a homotopy invariant functional which is an
analog of the notion of dimension of finite-dimensional vector spaces.

\begin{theorem}
{\em \cite{SaSt1,SaSt2}} There exists a unique additive homotopy invariant
functional
\[
d:\widehat{{\rm {Even}}}\left( M^{odd}\right) \rightarrow {Z}\left[ \frac
12\right] ,\quad {\mbox {and}}{\quad }d:\widehat{{\rm {Odd}}}\left(
M^{ev}\right) \rightarrow {Z}\left[ \frac 12\right] ,
\]
with the properties:
\end{theorem}

\begin{enumerate}
\item  {\em (relative dimension) }
\[
d\left( \widehat{L}+L_0\right) -d\left( \widehat{L}\right) ={\rm \dim }L_0
\]
{\em \ for a pair $\widehat{L}+L_0,\widehat{L}$ of subspaces, which differ
by a finite-dimensional space $L_0$;}

\item  {\em (complement) }
\begin{equation}
d\left( \widehat{L}\right) +d\left( \widehat{L}^{\perp }\right) =0,
\label{asavin}
\end{equation}
{\em \ \quad here }$\widehat{L}^{\perp }${\em \ denotes the orthogonal
complement of} $\widehat{L}.$
\end{enumerate}

It turns out that the spectral Atiyah--Patodi--Singer $\eta $-invariant is
equal to the dimension functional of subspaces.

\begin{theorem}
{\em \cite{SaSt1,SaSt2}$\!$ }Let $\!\widehat{L}_{+}\!\left( A\right) \!\!$
be the nonnegative spectral subspace of an elliptic self-adjoint
differential operator $A$ of positive order. Then the following equality
holds
\begin{equation}
\eta \left( A\right) =d\left( \widehat{L}_{+}\left( A\right) \right) ,
\label{dash}
\end{equation}
provided the order of the operator and the dimension of the manifold have
opposite parities.
\end{theorem}

The equality (\ref{dash}) is also valid for {\em admissible\/}
pseudodifferential operators in the sense of \cite{Gil7}.

\vspace{2mm}

{\bf 2.} The dimension functional of pseudodifferential subspaces enters the
index formula for elliptic operators in subspaces.

Indeed, consider two pseudodifferential subspaces $\widehat{L}_{1,2}\subset
C^\infty \left(M,E_{1,2}\right) $ and a pseudodifferential operator
\[
D:C^\infty \left( M,E_1\right) \longrightarrow C^\infty \left( M,E_2\right)
.
\]
Suppose that the following inclusion is valid $D\widehat{L}_1\subset
\widehat{L}_2.$ Then the restriction
\begin{equation}
D:\widehat{L}_1\longrightarrow \widehat{L}_2  \label{ttwostar}
\end{equation}
is called an {\em operator acting in subspaces}. The restriction
\begin{equation}
\sigma \left( D\right) :L_1\longrightarrow L_2  \label{ssymbl}
\end{equation}
of the principal symbol $\sigma \left( D\right) $ to the symbols of
subspaces over the cospheres $S^{*}M$ is called the {\em symbol of operator
in subspaces}. It can be shown that the closure of operator (\ref{ttwostar})
with respect to Sobolev norm defines an operator with Fredholm property if
and only if the symbol (\ref{ssymbl}) is {\em elliptic}, i.e. it is a vector
bundle isomorphism.

The following index formula for elliptic operators in subspaces was obtained
in \cite{SaSt1,SaSt2}.

\begin{theorem}
Let
\begin{equation}
\widehat{L}_{1,2}\in \widehat{{\rm {Even}}}\left( M^{odd}\right) \quad {%
\mbox {or}}\quad \widehat{{\rm {Odd}}}\left( M^{ev}\right) .  \label{chet}
\end{equation}
Then the index of an elliptic operator $D:\widehat{L}_1\longrightarrow
\widehat{L}_2$ is equal to
\begin{equation}
{\rm {ind}}\left( D,\widehat{L}_1,\widehat{L}_2\right) =\frac 12{\rm {ind}}%
\widetilde{D}+d\left( \widehat{L}_1\right) -d\left( \widehat{L}_2\right) ,
\label{inde}
\end{equation}
where in the case of odd subspaces
\[
\widetilde{D}:C^\infty \left( M,E_1\right) \longrightarrow C^\infty \left(
M,E_2\right)
\]
is the usual\footnote{%
That is, this operator acts in spaces of vector bundle sections.} elliptic
operator with principal symbol $\sigma (\widetilde{D})$ equal to
\[
\sigma \left( \widetilde{D}\right) =\sigma \left( D\right) \oplus \alpha
^{*}\sigma \left( D\right) :L_1\oplus \alpha ^{*}L_1\longrightarrow
L_2\oplus \alpha ^{*}L_2.
\]
For even subspaces, the following formula is valid
\[
\widetilde{D}:C^\infty \left( M,E_1\right) \longrightarrow C^\infty \left(
M,E_1\right) ,
\]
where the usual elliptic operator with principal symbol $\sigma (\widetilde{D%
})$ is defined by
\[
\sigma \left( \widetilde{D}\right) =\left[ \alpha ^{*}\sigma \left( D\right)
\right] ^{-1}\sigma \left( D\right) \oplus 1:L_1\oplus L_1^{\perp
}\longrightarrow L_1\oplus L_1^{\perp }.
\]
\end{theorem}

{\bf 3.} As a direct consequence of the index formula (\ref{inde}), we have
the following corollary.

\begin{corollary}
The fractional part of twice the functional $d$ is determined by the
principal symbol of the subspace $\widehat{L}$ as an element of the group $%
K\left( S^{*}M\right) /K\left( M\right) .$
\end{corollary}

Let us apply the index formula (\ref{inde}) for operators in subspaces to
the computation of the fractional part of the invariant $d.$ To this end let
us recall the following property of even (odd) subspaces.

\begin{theorem}
\label{th4} The symbol of a subspace $\widehat{L}$ with parity conditions
{\em (\ref{chet})} defines a $2$-torsion element in the group $K\left(
S^{*}M\right) /K\left( M\right) .$ In other words, for some number $N$ and a
vector bundle $F\in {\rm {Vect}}\left( M\right) $ on the base there exists
an isomorphism
\begin{equation}
\sigma :2^NL\longrightarrow \pi ^{*}F,\qquad 2^NL={\underbrace{L\oplus
\ldots \oplus L}_{2^N{\mbox {{\em\footnotesize copies}}}}}.  \label{omega}
\end{equation}
\end{theorem}

The proof of this theorem for even subspaces is contained in \cite{Gil7},
for odd subspaces --- in \cite{SaSt2}.

Consider the corresponding elliptic operator (in subspaces)
\begin{equation}
\widehat{\sigma }:2^N\widehat{L}\longrightarrow C^\infty \left( M,F\right)
\label{kukuk}
\end{equation}
with symbol (\ref{omega}). By virtue of (\ref{asavin}) the space of vector
bundle sections has ``dimension'' zero: $d\left( C^\infty \left( M,F\right)
\right) =0.$ Therefore, the index formula (\ref{inde}) gives an equality
\begin{equation}
{\rm {ind}}\left( \widehat{\sigma },2^N\widehat{L},C^\infty \left(
M,F\right) \right) =\frac 12{\rm {ind}}{\widetilde{\widehat\sigma }}%
+2^Nd\left( \widehat{L}\right) .  \label{omegasht}
\end{equation}
Placing the operators $\widehat{\sigma }\ $ and $\widetilde{\widehat\sigma }$
under the index sign, we can obtain the following expression for fractional
part of twice the value of the dimension functional (cf. \cite{SaScS3})
\begin{equation}
\left\{ 2d\left( \widehat{L}\right) \right\} =\frac 1{2^N}{\rm mod} \,2^N %
\mbox{-} {\rm {ind}}\left[ \left( 1\pm \alpha ^{*}\right) \widehat{\sigma }%
\right] ,  \label{moda}
\end{equation}
where ${\rm mod} \,2^N$-${\rm {ind}}D\in {Z}_{2^N}$ denotes the index of a
Fredholm operator $D$ reduced modulo $2^N,$ while $\alpha ^{*}\widehat{%
\sigma }$ denotes operator in subspaces
\[
\widehat{\alpha^*\sigma }:2^N\widehat{\alpha^*L}\longrightarrow
C^\infty(M,F).
\]

Note that the index of operator (\ref{kukuk}) as well as the index in (\ref
{moda}) as a residue modulo $2^N$ is determined by the principal symbol of
the operator. In the next section we show how this index can be computed.

\section{Index modulo $n$ and Pontryagin duality in $K$-theory}

Let us recall the main results of index theory modulo $n$ from \cite{SaScS3}.

\vspace{2mm}

{\bf 1.} Let $n$ be a natural number. {\em Elliptic operator modulo} $n$ is
an elliptic operator of the form
\begin{equation}
D:\,n\widehat{L}_1\oplus C^\infty \left( M,E_1\right) \rightarrow n\widehat{L%
}_2\oplus C^\infty \left( M,F_1\right) ,  \label{lab8}
\end{equation}
where
\[
\widehat{L}_1\subset C^\infty \left( M,E\right) ,\widehat{L}_2\subset
C^\infty \left( M,F\right)
\]
are pseudodifferential subspaces. The direct sums
\begin{equation}
n\left( \widehat{L}_1\oplus C^\infty \left( M,E_1\right) \right) \stackrel{nD%
}{\longrightarrow }n\left( \widehat{L}_2\oplus C^\infty \left( M,F_1\right)
\right)  \label{two}
\end{equation}
are called {\em trivial} operators modulo $n$. The group of stable homotopy
classes of elliptic operators modulo trivial ones is denoted by ${\rm {Ell}}%
\left( M,{Z}_n\right) .$ In \cite{SaScS3} it is shown that this group
defines $K$-theory with coefficients in ${Z}_n.$

\begin{theorem}
There is an isomorphism of groups
\[
{\rm {Ell}}\left( M,{Z}_n\right) \stackrel{\chi _n}{\simeq }K_c\left( T^{*}M,%
{Z}_n\right) .
\]
\end{theorem}

Here $K_c$ denotes $K$-theory with compact supports.

Let us give an explicit formula for the isomorphism $\chi_n$. Recall that $K$%
-theory with coefficients ${\bf Z}_n$ is defined by means of the so-called
Moore space $M_n$. This is a topological space with $K$-groups equal to
\[
\widetilde{K}^0\left( M_n\right) ={Z}_n,\quad K^1\left( M_n\right) =0.
\]
For instance, as a Moore space we can take a two-dimensional complex
obtained from the disk $D^2$ by identifying the points on its boundary under
the action of the group ${Z}_n:$%
\[
M_n=\left. \left\{ \left. D^2\subset {C}\right| \;\left| z\right| \leq
1\right\} \right/ \left\{ e^{i\varphi }\sim e^{i\left( \varphi +\frac{2\pi k}%
n\right) }\right\} .
\]
Denote the vector bundle corresponding to the generator of the group $%
\widetilde{K}^0\left(M_n\right) ={Z}_n$ by $\gamma _n$. Let us also fix a
trivialization $\beta$ of the direct sum $n\gamma _n$
\[
\begin{array}{ccc}
n\gamma _n & \stackrel{\beta }{\longrightarrow } & {C}^n.
\end{array}
\]
For a topological space $X$, its $K$-groups with coefficients ${Z}_n$ are
defined according to the formula
\begin{equation}
K^{*}\left( X,{Z}_n\right) =K^{*}\left( X\times M_n,X\times pt\right) .
\label{fdef}
\end{equation}

It is shown in \cite{SaScS3} that an arbitrary operator (\ref{lab8}) is
stably homotopic to some operator of the form
\begin{equation}
n\widehat{L}\stackrel{D}{\longrightarrow }C^\infty \left( M,F\right).
\label{simp1}
\end{equation}
An operator of this type defines a family of usual elliptic symbols on $M$
(here we use the difference construction for elliptic families)
\begin{equation}
\begin{array}{l}
\chi _n\left[ D\right] = \\
\left[ \pi ^{*}F\stackrel{\sigma ^{-1}\left( D\right) }{\longrightarrow }nL%
\stackrel{\beta ^{-1}\otimes 1_L}{\longrightarrow }\gamma _n\otimes nL%
\stackrel{1_\gamma \otimes \sigma \left( D\right) }{\longrightarrow }\gamma
_n\otimes \pi ^{*}F\right] \in K_c\left( T^{*}M\times M_n,T^{*}M\times
pt\right)
\end{array}
\label{long}
\end{equation}
The family is parametrized by the Moore space $M_n$.

Note that the index of an elliptic operator modulo $n$ is determined by its
principal symbol as a residue
\[
{\rm mod}\,n\mbox{-}{\rm {ind}}\left( n\widehat{L}\stackrel{D}{%
\longrightarrow }C^\infty \left( M,F\right) \right) \in {Z}_n.
\]
This index-residue can be computed topologically, the corresponding index
theorem is equivalent to the commutativity of the triangle
\begin{equation}
\begin{array}{ccc}
& \!\!\!\!\!{\rm {Ell}}\left( M,{Z}_n\right)  &  \\
\qquad \qquad \qquad \qquad \swarrow \mbox{\footnotesize $\chi_n$} &  &
\searrow \mbox{\footnotesize {\rm mod}\,$n$-{\rm ind}} \\
\qquad \qquad K_c\left( T^{*}M,{Z}_n\right)  & \stackrel{p_{!}}{%
-\!\!\!-\!\!\!-\!\!\!\longrightarrow } & \!\!\!\!\!\!\!\!\!{Z}_n,
\end{array}
\label{modn}
\end{equation}
where $p_{!}:K_c\left( T^{*}M,{Z}_n\right) \rightarrow {Z}_n$ is the direct
image mapping in $K$-theory with coefficients ${\bf Z}_n$.

We will use later an exact sequence relating elliptic operators modulo $n$
to the usual elliptic operators.

To this end let us denote the group of stable homotopy classes of elliptic
operators by ${\rm {Ell}}\left( M\right) $ (see \cite{AtSi1}), and a similar
stable homotopy group for pseudodifferential subspaces by ${\rm {Ell}}%
{}_1\left( M\right) $ (see \cite{SaScS3}).

The mappings
\[
{\rm {Ell}}\left( M\right) \stackrel{\chi _0}{\rightarrow }K_c\left(
T^{*}M\right) ,\quad {\rm {Ell}}{}_1\left( M\right) \stackrel{\chi _1}{%
\rightarrow }K_c^1\left( T^{*}M\right) ,
\]
which associate principal symbols to operators and subspaces, define
isomorphisms with the corresponding $K$-groups (recall (see \cite{APS0} or
\cite{SaScS3}) that the second mapping is defined as the composition
\[
\chi _1:\,{\rm {Ell}}{}_1\left( M\right) \rightarrow K\left( S^{*}M\right)
/K\left( M\right) \stackrel{\delta }{\rightarrow }K_c^1\left( T^{*}M\right) ,
\]
where the first mapping associates symbol to a subspace. Then we apply the
isomorphism $\delta $ induced from the coboundary mapping in $K$-theory:
\[
\delta :\,K\left( S^{*}M\right) \rightarrow K_c^1\left( T^{*}M\right) ).
\]
The following diagram is commutative
\begin{equation}
\begin{array}{ccccccccc}
{\rm {Ell}}\left( M\right)  & \stackrel{\times n}{\rightarrow } & {\rm {Ell}}%
\left( M\right)  & \stackrel{i}{\rightarrow } & {\rm {Ell}}\left( M,{Z}%
_n\right)  & \stackrel{\partial }{\rightarrow } & {\rm {Ell}}{}_1\left(
M\right)  & \stackrel{\times n}{\rightarrow } & {\rm {Ell}}{}_1\left(
M\right)  \\
\downarrow \mbox{\footnotesize $\chi _0$} &  & \downarrow %
\mbox{\footnotesize $\chi _0$} &  & \downarrow \mbox{\footnotesize $\chi _n$}
&  & \downarrow \mbox{\footnotesize $\chi _1$} &  & \downarrow %
\mbox{\footnotesize $\chi _1$} \\
K_c\left( T^{*}M\right)  & \stackrel{\times n}{\rightarrow } & K_c\left(
T^{*}M\right)  & \stackrel{i^{\prime }}{\rightarrow } & K_c\left( T^{*}M,{Z}%
_n\right)  & \stackrel{\partial ^{\prime }}{\rightarrow } & K_c^1\left(
T^{*}M\right)  & \stackrel{\times n}{\rightarrow } & K_c^1\left(
T^{*}M\right) ,
\end{array}
\label{vlong}
\end{equation}
here the mapping $\times n$ sends an element $x$ of an abelian group to its
multiple $nx,$ the mapping $i$ is induced by the embedding of usual elliptic
operators into the set of elliptic operators modulo $n,$ while the Bockstein
homomorphism $\partial $ is defined by the formula
\[
\partial \left[ n\widehat{L}_1\oplus C^\infty \left( M,E_1\right) \stackrel{D%
}{\longrightarrow }n\widehat{L}_2\oplus C^\infty \left( M,F_1\right) \right]
=\left[ \widehat{L}_1\right] -\left[ \widehat{L}_2\right] .
\]

Consider the natural inclusions ${Z}_n\subset {Z}_{mn}$ for a pair of
natural numbers $n,m$. The direct limit
\[
\mathop{{\rm lim}\,{\bf Z}_n} \limits_{\stackrel{\longrightarrow}{%
n\rightarrow\infty}}= {Q}/{Z}
\]
permits us to define $K$- and ${\rm {Ell}}$-groups with coefficients in ${Q}/%
{Z}$:
\begin{eqnarray*}
K_c\left( T^{*}M,{Q}/{Z}\right) &= &\mathop{{\lim}\,K_c}\limits_{%
\longrightarrow}\left( T^{*}M,{Z}_n\right) , \\
{\rm {Ell}}\left( M,{Q}/{Z}\right) &= &\mathop{{\lim}\,{\rm Ell}}%
\limits_{\longrightarrow}\left( M,{Z}_n\right) .
\end{eqnarray*}
Moreover, (\ref{vlong}) transforms into the diagram with exact rows
\begin{equation}
\begin{array}{ccccccccc}
\!{\rm {Ell}}\left( M\right) & \stackrel{}{\!\!\!\rightarrow\!\!\! } & {\rm {%
Ell}}\left( M\right) \otimes {Q} & \stackrel{i}{\!\!\!\rightarrow\!\!\! } &
{\rm {Ell}}\left( M,{Q}/{Z}\right) & \stackrel{\partial}{\!\!\!\rightarrow\!%
\!\! } & {\rm {Ell}}{}_1\left( M\right) & \stackrel{}{\!\!\!\rightarrow\!\!%
\! } & {\rm {Ell}}{}_1\left( M\right) \otimes {Q} \! \\
\;\downarrow \mbox{\footnotesize $\chi _0$} &  & \qquad \downarrow %
\mbox{\footnotesize $\chi _0\otimes 1$} &  & \quad\downarrow %
\mbox{\footnotesize $\chi _{{\mathbf{}Q}/{\mathbf{}Z}}$} &  & \downarrow %
\mbox{\footnotesize $\chi _1$} &  & \qquad\downarrow \mbox{\footnotesize
$\chi _1\otimes 1$} \\
\!K_c\left( T^{*}M\right) & \stackrel{}{\!\!\!\rightarrow\!\!\! } & K_c\left(
T^{*}M\right) \otimes {Q} & \stackrel{i^{\prime }}{\!\!\!\rightarrow\!\!\! }
& K_c\left( T^{*}M,{Q}/{Z}\right) & \stackrel{\partial^{\prime }}{%
\!\!\!\rightarrow\!\!\! } & K_c^1\left( T^{*}M\right) & \stackrel{}{%
\!\!\!\rightarrow\!\!\! } & K_c^1\left( T^{*}M\right) \otimes {Q}.\!
\end{array}
\label{longa}
\end{equation}

{\bf 2.} Consider the intersection form
\begin{equation}
K_c^i\left( T^{*}M,{Q/Z}\right) \times K_c^i\left( M\right) \longrightarrow
K_c^0\left( T^{*}M,{Q/Z}\right) \stackrel{p_{!}}{\rightarrow }{Q/Z,}
\label{aa}
\end{equation}
which is induced by the product and the direct image mapping $p_{!}$
in $K$-theory. The intersection of elements $x$ and $y$ is denoted by $x\cap
y.$

Recall that a pairing
\[
\left\langle \cdot,\cdot \right\rangle :G_1\times G_2\rightarrow G_3
\]
of abelian groups $G_{1,2}$ with values in abelian group $G_3$ is called
{\em nondegenerate} if the mappings
\[
\left\langle x,\cdot \right\rangle :G_2\rightarrow G_3{\mbox { and }}%%
\left\langle \cdot ,y\right\rangle :G_1\rightarrow G_3,
\]
obtained by fixing the values of one of the arguments are zero only for $x=0$
(correspondingly $y=0$).

\begin{theorem}
\label{th111}{\em (Pontryagin duality)} The pairing {\em (\ref{aa})} is
nondegenerate. In addition, fixing its first argument, we obtain an
isomorphism
\begin{equation}
K_c^i\left( T^{*}M,{Q/Z}\right) \simeq {\rm {Hom}}\left( K^i\left( M\right) ,%
{Q/Z}\right) {.}  \label{iso4}
\end{equation}
\end{theorem}

\noindent {\em Proof. } It can be easily shown that the isomorphism (\ref
{iso4}) implies the nondegeneracy of the pairing.

Similarly to Poincar\'e duality (e.g., see \cite{BoTu1}), the Pontryagin
duality can be proved by means of the Mayer--Vietoris principle. For an
arbitrary open subset $U\subset M$, let us consider the mapping
\begin{equation}
K_c^i\left( T^{*}U,{Q/Z}\right) \longrightarrow {\rm {Hom}}\left( K^i\left(
U\right) ,{Q/Z}\right) .  \label{sueta}
\end{equation}
We want to prove that this mapping is an isomorphism for $U=M$. According to
the Mayer--Vietoris principle it suffices a) to verify it for a contractible
subset $U;$ b) to prove the inductive statement: if the mapping is an
isomorphism for two open subsets $U,V$ and for their intersection $U\cap V$
then it is an isomorphism for the union.

Consider a contractible set $U.$ Then
\[
K_c^{*}\left( T^{*}U,{Q/Z}\right) =K^{*}\left( pt\right) ={Q/Z}\oplus 0
\]
and
\[
K^{*}\left( U\right) ={Z}\oplus 0,
\]
while the corresponding mapping
\[
{Q/Z\rightarrow }{\rm {Hom}}\left( {Z},{Q/Z}\right)
\]
is an isomorphism. The validity of condition a) is thereby proved.

To prove b) consider the diagram
\[
\begin{array}{ccccccc}
\!\!\to\!\!\!\!\!\!\!\! & K_c^{i-1}\left( T^{*}\left( U\!\bigcap\! V\right),{%
Q/Z} \right) & \!\!\!\!\!\!\!\!\to\!\!\!\!\!\!\!\! & K_c^{i}\left(
T^{*}\left( U\!\bigcup\! V\right),{Q/Z} \right) & \!\!\!\!\!\!\!\!\to\!\!\!%
\!\!\!\!\! & K_c^{i}\left( T^{*}U\!\sqcup\! T^{*}V,{Q/Z}\right) &
\!\!\!\!\!\!\!\!\to\!\! \\[2mm]
& \downarrow &  & \downarrow &  & \downarrow &  \\[2mm]
\!\!\to\!\!\!\!\! & {\rm {Hom}}\left( K^{i-1}\!\left( U\!\bigcap\! V\right) ,{%
Q/Z}\right) & \!\!\!\!\!\!\!\!\to\!\!\!\!\!\!\!\! & {\rm {Hom}}\left(
K^{i}\!\left( U\!\bigcup\! V\right) ,{Q/Z}\right) & \!\!\!\!\!\!\!\!\to\!\!\!%
\!\!\!\!\! & {\rm {Hom}}\left( K^{i}\!\left( U\!\sqcup\! V\right) ,{Q/Z}\right)
& \!\!\!\!\to\!.
\end{array}
\]
Here the upper row is the Mayer--Vietoris exact sequence, the lower row is
obtained from a similar sequence applying the functor ${\rm {Hom}}%
\left(\cdot,{Q/Z}\right) $ (this functor preserves the exactness of
sequences). Vertical mappings are induced by the intersection pairing. This
diagram is commutative up to sign. Suppose that the mapping (\ref{sueta}) is
an isomorphism for $U,V$ and their intersection $U\cap V$. In accordance
with the 5-lemma we obtain the isomorphism for the union $U\cup V$. Thus, we
establish condition b).

Therefore, both conditions of the Mayer--Vietoris principle are satisfied.
The theorem is proved.

\begin{remark}
{\em For }$K${\em -theory with coefficients in a topological group }${R/Z}$%
{\em \ (its definition can be found in \cite{APS0}), one can obtain by the
same method the Pontryagin duality of groups }$K_c^i\left( T^{*}M,{R/Z}%
\right) ${\em \ and } $K^i\left( M\right) $:
\[
K_c^i\left( T^{*}M,{R/Z}\right) \simeq {\rm {Hom}}\left( K^i\left( M\right) ,%
{R/Z}\right) ,{\mbox { }}K^i\left( M\right) \simeq {\rm {Hom}}\left(
K_c^i\left( T^{*}M,{R/Z}\right) ,{R/Z}\right) ,
\]
{\em \ i.e. both groups are character groups of each other.}
\end{remark}

\begin{remark}
{\em By a similar method one can prove the Pontryagin duality for an
arbitrary }$K${\em -oriented closed manifold or manifold with boundary
(i.e., a manifold with a }$spin^c${\em \ structure in the tangent bundle).
From this more general viewpoint, Theorem \ref{th111} establishes  the
Pontryagin--Lefschetz duality in }$K${\em -theory} {\em \ for an
almost-complex manifold }$T^{*}M${\em \ of groups with compact supports and
absolute groups.}
\end{remark}

\begin{remark}
\label{zama1}{\em One could also prove the nondegeneracy of the pairing}
\[
K_c^i\left( T^{*}M\right) \times K^i\left( M,{Q/Z}\right) \longrightarrow
K_c^0\left( T^{*}M,{Q/Z}\right) \stackrel{p_{!}}{\rightarrow }{Q/Z.}
\]
\end{remark}

Just like in (co)homology theory (see \cite{Pon2}), the Pontryagin duality
implies Poincar\'e duality for torsion subgroups. The bilinear form
defining this duality is called the {\em linking form}.

\section{Linking form in $K$-theory}

\label{parg} \noindent
{\bf 1.} Consider the Bockstein homomorphism
\[
\begin{array}{cc}
\partial :K_c^i\left( T^{*}M,{Q}/{Z}\right)  & \rightarrow K_c^{i-1}\left(
T^{*}M\right)
\end{array}
\]
(see diagram (\ref{longa})). The range of this mapping consists of finite
order elements. Denote by ${\rm {Tor}}G$ the torsion subgroup of an abelian
group $G.$

\begin{definition}
The linking form {\em is the pairing}
\begin{eqnarray*}
\cap :{\rm {Tor}}K_c^{i-1}\left( T^{*}M\right) \times {\rm {Tor}}K^i\left(
M\right) \longrightarrow {Q}/{Z,} \\
(x,y)\mapsto x^{\prime }\cap y,
\end{eqnarray*}
{\em where }$x^{\prime }\in K_c^i\left( T^{*}M,{Q}/{Z}\right) ${\em \ is an
arbitrary element of }$K${\em -theory with coefficients such that }$\partial
x^{\prime }=x,${\em \ where } $x^{\prime }\cap y\in {Q}/{Z}$%%
{\em \ is the intersection index from the previous section.}
\end{definition}

Similarly to the intersection form, the linking form is defined by a product
\begin{eqnarray}
{\rm {Tor}}K_c^{i-1}\left( T^{*}M\right) \times {\rm {Tor}}K^i\left(
M\right) & \longrightarrow & K_c\left( T^{*}M,{Q}/{Z}\right) {,}  \nonumber
\\
(x,y) & \mapsto & x^{\prime }y,\;\partial x^{\prime }=x  \label{mumu}
\end{eqnarray}

\begin{lemma}
The product {\em (\ref{mumu}) } and the linking form are well-defined.
\end{lemma}

\noindent {\em Proof. }We need to show that the indeterminacy of the choice
of $x^{\prime }$ does not affect the product $x^{\prime }y.$ Let $x=\partial
x^{\prime \prime }.$ For the difference $x^{\prime }-x^{\prime \prime }$ we
obtain $\partial \left( x^{\prime }-x^{\prime \prime }\right) =0.$ Thus, $%
x^{\prime }-x^{\prime \prime }=i\left( z\right) $ for some $z\in
K_c^i\left(T^{*}M\right) \otimes {Q}.$ However, the product $zy\in
K_c^0\left( T^{*}M\right) \otimes {Q}$ is a torsion element. Hence, $zy=0.$
We have proved that $x^{\prime }y=x^{\prime \prime }y.$ Lemma is proved.

The linking form can be defined similarly using the torsion property of the
second argument. To this end consider the Bockstein homomorphism $\partial $
in the coefficient sequence
\[
\begin{array}{cccccc}
\ldots \rightarrow  & K^{i-1}\left( M,{Q/Z}\right)  & \stackrel{\partial }{%
\rightarrow } & K^i\left( M\right)  & \rightarrow  & K^i\left( M\right)
\otimes {Q}\rightarrow \ldots
\end{array}
.
\]
As a linking index of elements $x\in {\rm {Tor}}K_c^{i-1}\left(
T^{*}M\right) ,y\in {\rm {Tor}}K^i\left( M\right) $ let us put
\[
x\cap ^{\prime }y=x\cap y^{\prime }.
\]
It turns out that both methods define linking pairings which coincide up to
sign. Namely, the following equality is valid
\[
x\cap ^{\prime }y=\left( -1\right) ^{\deg x+1}x\cap y.
\]
This formula follows from the proposition.

\begin{proposition}
Consider $x^{\prime }\in K_c^i\left( T^{*}M,{Q/Z}\right) ,y^{\prime }\in
K^j\left( M,{Q/Z}\right) $. Then the following equality holds
\[
\partial x^{\prime }y^{\prime }=\left( -1\right) ^{\deg x^{\prime
}+1}x^{\prime }\partial y^{\prime }.
\]
\end{proposition}

\noindent {\em Proof. }Let us assume that $x^{\prime },y^{\prime }$ are
induced from $K$-groups with coefficients ${Z}_n$ for $n$ large enough:
\begin{equation}
x^{\prime }=I_{*}x_0,y^{\prime }=n_{*}y_0,\;x_0\in K_c^i\left( T^{*}M,{Z}%
_n\right) ,y_0\in K^j\left( M,{Z}_n\right) ,  \label{ax}
\end{equation}
where $I$ denotes the embedding ${Z}_n\subset {Q/Z.}$

Consider the exact sequence
\[
0\rightarrow {Z}_n\stackrel{\times n}{\longrightarrow }{Z}_{n^2}\stackrel{%
{\rm mod}\,n}{\longrightarrow }{Z}_n\rightarrow 0
\]
and the corresponding sequence in $K$-theory
\[
\ldots \rightarrow K_c^{i+j}\left( T^{*}M,{Z}_n\right) \stackrel{\partial
^{\prime \prime }}{\longrightarrow }K_c^{i+j+1}\left( T^{*}M,{Z}_n\right)
\stackrel{\times n}{\longrightarrow }K_c^{i+j+1}\left( T^{*}M,{Z}%
_{n^2}\right) \rightarrow \ldots .
\]
%Let us act on the product
%$x_0y_0\in K_c^{i+j}\left( T^{*}M,{\mathbf{}Z}_n\right)$
The Bockstein homomorphism satisfies the Leibniz rule
\begin{equation}
\partial ^{\prime \prime }\left( x_0y_0\right) =\partial x_0y_0+\left(
-1\right) ^{\deg x_0}x_0\partial y_0.  \label{ay}
\end{equation}
On the other hand, by virtue of the exactness, we have $\times n\circ
\partial ^{\prime \prime }\left( x_0y_0\right) =0.$ Hence,
\begin{equation}
I_{*}\partial ^{\prime \prime }\left( x_0y_0\right) =0.  \label{az}
\end{equation}
The expressions (\ref{ax}), (\ref{ay}), (\ref{az}) together imply the
desired
\[
\partial x^{\prime }y^{\prime }+\left( -1\right) ^{\deg x^{\prime
}}x^{\prime }\partial y^{\prime }=0.
\]

\begin{theorem}
\label{mylo} {\em (Poincar\'e duality for torsion groups)} The linking form
\[
\begin{array}{ccccc}
{\rm {Tor}}K_c^{i-1}\left( T^{*}M\right)  & \times  & {\rm {Tor}}K^i\left(
M\right)  & \longrightarrow  & {Q}/{Z}
\end{array}
\]
is nondegenerate. In particular, it defines isomorphisms
\begin{eqnarray*}
{\rm {Tor}}K_c^{i-1}\!\left( T^{*}\!M\right) \!\simeq {\rm {Hom}}\!\left(
{\rm {Tor}}K^i\!\left( M\right) \!,{Q}/{Z}\right) , \\[2mm]
{\rm {Tor}}K^i\!\left( M\right) \simeq {\rm {Hom}}\!\left( {\rm {Tor}}%
K_c^{i-1}\!\left( T^{*}\!M\right) \!,{Q}/{Z}\right) .
\end{eqnarray*}
\end{theorem}

\begin{corollary}
The groups ${\rm {Tor}}K_c^{i-1}\left( T^{*}M\right) $ and ${\rm {Tor}}%
K^i\left( M\right) $ are isomorphic.
\end{corollary}

This follows from (noncanonical) isomorphism $G\simeq {\rm {Hom}}\left( G,{Q}%
/{Z}\right) $ for a finite abelian group $G.$

\noindent {\em Proof of Theorem} \ref{mylo}. Let us first prove the
nondegeneracy with respect to the second argument. Suppose that $x\cap y=0$
for an arbitrary $x\in {\rm {Tor}}K_c^{i-1}\left( T^{*}M\right) .$ Therefore,
for any $x^{\prime }\in K_c^i\left( T^{*}M,{Q}/{Z}\right) $ we
also have $x^{\prime }\cap y=0.$ The Pontryagin duality implies that $y=0.$

The nondegeneracy of the pairing with respect to the first argument follows
from the Pontryagin duality corresponding to the intersection form
\[
K_c^{i-1}\left( T^{*}M\right) \times K^i\left( M,{Q}/{Z}\right) \rightarrow {%
Q}/{Z}
\]
(see Remark \ref{zama1}).

The second claim of the theorem follows from the finiteness of the torsion
subgroups. The theorem is thereby proved.

\vspace{2mm}

{\bf 2.} The isomorphism of elliptic theory and $K$-theory from the first
part of the previous section permits us to define the linking form in terms
of elliptic operators.

\begin{definition}
The linking form in elliptic theory {\em is the bilinear pairing }
\begin{eqnarray*}
{\rm {Tor}}{\rm {Ell}}_1\left( M\right) \times {\rm {Tor}}K^0\left( M\right)
\longrightarrow {Q}/{Z,} \\
(x,y)\mapsto \qquad {\rm {ind}}x^{\prime }y,\;\partial x^{\prime }=x.
\end{eqnarray*}
\end{definition}

\begin{proposition}
The linking forms in elliptic and $K$-theory are isomorphic, i.e. the
following diagram commutes
\[
\begin{array}{ccc}
{\rm {Tor}}{\rm {Ell}}_1\left( M\right) \times {\rm {Tor}}K^0\left( M\right)
& \rightarrow  & {Q}/{Z} \\
\downarrow \mbox{\footnotesize $\chi_1\times 1$} &  & \downarrow %
\mbox{\footnotesize\rm Id} \\
{\rm {Tor}}K_c^1\left( T^{*}M\right) \times {\rm {Tor}}K^0\left( M\right)  &
\rightarrow  & {Q}/{Z}.
\end{array}
\]
\end{proposition}

\noindent {\em Proof} of the proposition follows from the commutative
diagram
\[
\begin{array}{ccccc}
{\rm {Tor}}{\rm {Ell}}_1\left( M\right) \times {\rm {Tor}}K^0\left( M\right)
& \rightarrow & {\rm {Ell}}\left( M,{Q}/{Z}\right) & \stackrel{{\rm {ind}}}{%
\rightarrow } & {Q}/{Z} \\
\downarrow &  & \downarrow &  & \parallel \\
{\rm {Tor}}K_c^1\left( T^{*}M\right) \times {\rm {Tor}}K^0\left( M\right) &
\rightarrow & K_c\left( T^{*}M,{Q}/{Z}\right) & \stackrel{p_{!}}{\rightarrow
} & {Q}/{Z.}
o\end{array}
\]
The left square commutes due to the isomorphism of coefficient sequences in $%
K$-theory and in elliptic theory (see (\ref{longa})), while the right square
commutes by virtue of the modulo $n$ index formula (see (\ref{modn})). The
proposition is proved.

Let us write an explicit formula for the linking form in elliptic theory.
Let $x=\left[ \widehat{L}\right]$ be a pseudodifferential subspace and $%
y=\left[ G_1\right] -\left[ G_2\right]$ be a difference of vector bundles $%
G_{1,2}\in {\rm {Vect}}\left( M\right) .$ According to the first definition
of the linking form for some element $x^{\prime }=\left[ n\widehat{L}%
\stackrel{\widehat{\sigma }}{\rightarrow }C^\infty \left( M,F\right) \right]
\in {\rm {Ell}}\left( M,{Z}_n\right) $ the linking coefficient is equal to
\begin{equation}
x\cap y={\rm {ind}}x^{\prime }y=\frac 1n{\rm mod} \,n\mbox{-}{\rm {ind}}%
\left( \widehat{\sigma }\otimes 1_{G_1-G_2}\right) ,  \label{ona}
\end{equation}
where $\widehat{\sigma }\otimes 1_{G_1-G_2}$ denotes operator $\widehat{%
\sigma }$ with coefficients in $G_1-G_2.$ By the second definition of the
linking form for some $y^{\prime }=\left[ mG_1\stackrel{g}{\rightarrow }%
mG_2\right] \in K^1\left( M,{Z}_m\right) $ (see \cite{APS0}) this index is
defined as
\[
x\cap y={\rm {ind}}xy^{\prime }=\frac 1m{\rm mod} \,m\mbox{-}{\rm {ind}}%
\left( 1_{\widehat{L}}\otimes g\right),
\]
where the operator in subspaces $1_{\widehat{L}}\otimes g$ has principal
symbol $mL\otimes G_1\stackrel{1_L\otimes g}{\longrightarrow} mL\otimes G_2$.

\begin{remark}
{\em Applying the Poincar\'e isomorphism in complex }$K${\em -theory to the
manifold }$T^{*}M${\em \ (see \cite{BaDo2}) }
\[
K_c^1\left( T^{*}M\right) \simeq K_1\left( M\right) ,
\]
{\em the linking pairing can be considered as a }nondegenerate{\em \ pairing
of }$K${\em -homology with} $K${\em -cohomology of the manifold }$M${\em : }
\[
\begin{array}{ccc}
{\rm {Tor}}K_1\left( M\right) \times {\rm {Tor}}K^0\left( M\right)  &
\stackrel{\cap }{\rightarrow } & {Q}/{Z.}
\end{array}
\]
\end{remark}

Let us also note the similarity of expression (\ref{ona}) for the linking
index and (\ref{moda}) for the fractional part of the dimension functional.
Unlike the linking index, the latter formula contains the action of the
antipodal involution instead of the product with some bundle. We show in the
following section that the antipodal involution acts on $K$-groups as a
tensor product with the orientation bundle of the manifold (cf. \cite{Karo2}
in the orientable case).

\section{Antipodal involution and orientability}

Let $V$ be a real vector bundle over a compact space $X.$ Consider the
antipodal involution
\[
\begin{array}{rcc}
\alpha :V & \longrightarrow & V, \\
v & \mapsto & -v.
\end{array}
\]

\begin{theorem}
\label{thm101}The induced involution $\alpha ^{*}$ in $K$-theory is equal to
\begin{equation}
\alpha ^{*}=\left( -1\right) ^n\Lambda ^n\left( V\right) :K_c^{*}\left(
V\right) \longrightarrow K_c^{*}\left( V\right) ,  \label{star3}
\end{equation}
where $\Lambda ^n\left( V\right) $ for $n=\dim V$ is the orientation bundle
of $V.$
\end{theorem}

\noindent {\em Proof.} Without loss of generality we can assume that $V$ is
even-dimensional: $n=\dim V=2k.$ We can also consider only the action of $%
\alpha ^{*}$ on even $K$-groups $K_c^0\left( V\right) $ (the odd case
reduces to this situation taking a sum with a one-dimensional trivial
bundle, since in this case both sides of the equality (\ref{star3}) change
by $-1$).

M.~Karoubi in \cite{Karo1,Karo2} found a description of the groups $%
K_c^0\left( V\right) $ in terms of Clifford algebras. Let us recall the
basic definitions. On the space $X$, consider the bundle $Cl\left( V\right) $
of Clifford algebras of the vector bundle $V.$ Over a point $x$ of the base $%
X$ the Clifford algebra $Cl\left( V_x\right) $ is multiplicatively generated
by the vector space $V_x$ with the following relations
\[
v_1v_2+v_2v_1=2\left\langle v_1,v_2\right\rangle
\]
for some scalar product $\left\langle,\right\rangle $ in $V.$

Let us consider quadruples $\left(E,c,f_1,f_2\right) ,$ where $E$ is a
complex vector bundle over $X,$ $c$ is a homomorphism
\[
c:Cl\left( V\right) \longrightarrow {\rm {Hom}}\left( E,E\right)
\]
of algebra bundles. We say that $c$ defines a Clifford module structure on $E
$. The involutions $f_{1,2}$ of $E$ are supposed to skew commute with the
Clifford multiplication:
\[
f_ic(v)+c(v)f_i=0,\qquad i=1,2.
\]
Stable homotopy equivalence relation is defined on the set of
quadruples $\left( E,c,f_1,f_2\right)$. The corresponding group of
equivalence classes is denoted by $K^V\left(X\right) .$ It is proved in \cite
{Karo1} that this group is isomorphic to $K_c^0\left( V\right) .$ The
isomorphism
\[
t:K^V\left( X\right) \longrightarrow K_c^0\left( V\right)
\]
is given by the following explicit formula
\begin{equation}
t\left[ E,c,f_1,f_2\right] =\!\left[ \pi ^{*}\!\ker \left( f_1\!-\! 1\right)
\stackrel{\left( 1-c\left( v\right) f_2\right) \left( 1+c\left( v\right)
f_1\right) }{ -\!\!\!-\!\!\!-\!\!\!-\!\!\! \longrightarrow }\pi ^{*}\!\ker
\left( f_2\!-\! 1\right) \right] ,\quad\! {\mbox {where}}\!\quad \pi :SV
\!\rightarrow\! X.
\label{tt}
\end{equation}
Here $SV$ is the sphere bundle of $V,$ and the element in the left hand side
of the equality is understood in the sense of difference construction (e.g.,
see \cite{Mis1}) for the relative group $K\left( BV,SV\right) \simeq
K_c\left(V\right) .$

The formula (\ref{tt}) implies that the antipodal involution $\alpha ^{*}$
acts on the quadruples according to the formula
\[
\alpha ^{*}\left[ E,c,f_1,f_2\right] =\left[ E,-c,f_1,f_2\right] .
\]
Let us show that $\left[ E,-c,f_1,f_2\right] $ differs from the quadruple $%
\left[ E,c,f_1,f_2\right] \otimes \Lambda ^{2k}\left( V\right) $ by a vector
bundle isomorphism (cf. \cite{Gil8}).

To this end, consider the local orthonormal base $e_1,e_2,...,e_{2k}$ of $V$%
. Let us define the element
\[
\beta =i^kc\left( e_1\right) \ldots c\left( e_{2k}\right) .
\]
One verifies easily that $\beta^2=1$ and $\beta $ skew commutes with the
Clifford action $c$ and commutes with the involutions $f_{1,2}$%
\begin{equation}
\beta c(v)+c(v)\beta =0,\quad f_i\beta =\beta f_i.  \label{ssttar}
\end{equation}
A direct computation shows that the change of the orthonormal base changes $%
\beta $ by the sign of the transition matrix determinant. Therefore,
globally this element defines a vector bundle isomorphism
\[
\beta :E\longrightarrow E\otimes \Lambda ^{2k}\left( V\right)
\]
(for brevity we denote it by the same symbol). At the same time the
commutation relations (\ref{ssttar}) transform to
\[
\beta ^{-1}\left( c(v)\otimes 1_{\Lambda ^{2k}\left( V\right) }\right) \beta
=-c(v),\quad \beta ^{-1}\left( f_i\otimes 1_{\Lambda ^{2k}\left( V\right)
}\right) \beta =f_i.
\]
Hence, the quadruples
\[
\alpha ^{*}\left[ E,c,f_1,f_2\right] \quad {\mbox {and}}\quad \left[
E,c,f_1,f_2\right] \otimes \Lambda ^{2k}\left( V\right)
\]
are isomorphic. Theorem is thereby proved.

\begin{corollary}
The same formula holds for $K$-theory with coefficients ${\bf Z}_n$
\[
\alpha ^{*}=\left( -1\right) ^n\Lambda ^n\left( V\right) :K_c^{*}\left( V,{Z}%
_n\right) \longrightarrow K_c^{*}\left( V,{Z}_n\right) .
\]
\end{corollary}

Indeed, $K$-theory with coefficients is defined by means of the Moore
space $M_n$ by the formula
\[
K_c^{*}\left( V,{Z}_n\right) =K_c^{*}\left( V\times M_n,V\times pt\right) .
\]
Let us apply our theorem to the space $X\times M_n$ and the pull-back of
bundle $V$. This gives the desired formula for $K$-theory with coefficients.

\begin{remark}
{\em In the case when }$X${\em \ is a smooth manifold with the cotangent
bundle }$V=T^{*}X,${\em \ (\ref{tt}) defines a class of elliptic symbols,
such that an arbitrary symbol reduces to a symbol of this form by a stable
homotopy.}
\end{remark}

\section{The main theorem}

\begin{theorem}
Fractional part $\{2d(\hat{L})\}$ of twice the value of the dimension
functional $d$ on subspace $\widehat{L}$ is equal to the linking index of
the subspace with the orientation bundle of the manifold $\Lambda ^n\left(
M\right) ,n=\dim M$:
\[
\left\{ 2d\left( \widehat{L}\right) \right\} =\left[ L\right] \cap \left(
1-\left[ \Lambda ^n\left( M\right) \right] \right) \in {Z}\left[ \frac
12\right] /{Z.}
\]
\end{theorem}

\noindent {\em Proof.} Recall the expression (\ref{moda}) of the fractional
part of the dimension functional $d$:
\[
\left\{ 2d\left( \widehat{L}\right) \right\} =\frac 1{2^N}{\rm mod} \,2^N %
\mbox{-}{\rm {ind}}\left[ \left( 1\pm \alpha ^{*}\right)\hat \sigma \right]
,\quad \left[ \sigma \right] \in K_c^0\left( T^{*}M,{Z}_{2^N}\right) ,
\]
where $\sigma :2^NL\rightarrow \pi ^{*}F$ is a vector bundle isomorphism. By
Theorem \ref{thm101} of the previous section we have
\[
\left( 1\pm \alpha ^{*}\right) \left[ \sigma \right] =\left[ 1-\Lambda
^n\left( M\right) \right] \left[ \sigma \right] .
\]
Hence,
\[
\left\{ 2d\left( \widehat{L}\right) \right\} =\frac 1{2^N}{\rm mod} \,2^N %
\mbox{-}{\rm {ind}}\left[ \left( 1-\Lambda ^n\left( M\right) \right)
\hat\sigma \right].
\]
This coincides with the definition of the linking index in Section \ref{parg}%
, see (\ref{ona}).

Note also that $[L]\in K_c^1(T^*M)$ is a torsion element by virtue of
Theorem \ref{th4}, the torsion property for the difference $[1-\Lambda
^n\left( M\right)]\in K(M)$ is proved in Proposition \ref{thr3}. The theorem
is proved.

\begin{corollary}
The dimension functional takes half-integer values
\[
\left\{ 2d\left( \widehat{L}\right) \right\} =0
\]
when the manifold is orientable.
\end{corollary}

\begin{proposition}
\label{thr3} Let $M$ be a nonorientable manifold of dimension $2k$ or $2k+1$%
. Then the following integrality estimate of the invariant $d$ is
valid:
\begin{equation}
\left\{ 2^{k+1}d\left( \widehat{L}\right) \right\} =0.  \label{eval}
\end{equation}
\end{proposition}

The same statements hold for the $\eta$-invariant.
%It is shown in the next section that the above
%estimate is exact for $k\geq 1$.

\noindent {\em Proof of Proposition} \ref{thr3}. The orientation bundle $%
\Lambda ^n\left( M^n\right) $ is a one-dimensional bundle with structure
group ${Z}_2.$ Hence, it is the pull-back of the universal bundle from the
classifying space ${BZ}_2={RP}^\infty ,$ i.e. there is a vector bundle
isomorphism
\[
\Lambda ^n\left( M^n\right) \simeq f^{*}\gamma ,\quad {%
\mbox {for some
mapping }}f:M^n\rightarrow {RP}^N.
\]
Here $\gamma $ is the line bundle over the projective space ${RP}^N.$ By the
approximation theorem we can suppose that $N=n.$ The reduced $K$-groups of
the projective spaces are well known
\[
\widetilde{K}\left( {RP}^{2k}\right) \simeq \widetilde{K}\left( {RP}%
^{2k+1}\right) \simeq {Z}_{2^k}.
\]
Thus, we have
\[
2^k\left( 1-\left[ \Lambda ^n\left( M^n\right) \right] \right) =0.
\]
Hence, we obtain the desired
\[
\left\{ 2^{k+1}d\left( \widehat{L}\right) \right\} =\left[ L\right]
\cap\left[ 2^k\left( 1-\Lambda ^n\left( M^n\right) \right) \right] =\left[
L\right] \cap 0=0.
\]

\section{Examples}

{\bf 1.} Consider the even-dimensional real projective space $\!{RP}^{2n}.$
The reduced $K$-group of this manifold is cyclic
\[
\widetilde{K}\left( {RP}^{2n}\right) \simeq {Z}_{2^n}.
\]
The generator is given by the orientation bundle
\[
1-\left[ \Lambda ^{2n}\left( {RP}^{2n}\right) \right] \in \widetilde{K}%
\left( {RP}^{2n}\right) .
\]
On the other hand, the projective space ${RP}^{2n}$ has a $pin^c$ structure,
while the principal symbol of the self-adjoint $pin^c$ Dirac operator $D$ on
it (this operator was constructed in \cite{Gil8}) is a generator of the
isomorphic group
\[
\left[ \sigma \left( D\right) \right] \in K_c^1\left( T^{*}{RP}^{2n}\right) =%
{\rm {Tor}}K_c^1\left( T^{*}{RP}^{2n}\right) \simeq {Z}_{2^N}.
\]
The nondegeneracy of the linking form implies that these generators have a
nontrivial linking index
\[
2^{n-1}\left[ \sigma \left( D\right) \right] \cap \left[ 1-\left[
\Lambda^{2n}\left( {RP}^{2n}\right) \right] \right] =\frac 12.
\]
Hence, the $\eta $-invariant of the $pin^c$ Dirac operator $D$ has a large
fractional part (see \cite{Gil8}):
\[
\left\{ 2^n\eta \left( D\right) \right\} = 2^{n-1}\left[ \sigma \left(
D\right) \right] \cap \left[ 1-\left[ \Lambda^{2n}\left( {RP}^{2n}\right)
\right] \right] =\frac 12.
\]
This example shows that the estimate (\ref{eval}) is precise
on even-dimensional manifolds.

\vspace{2mm}

{\bf 2.} {\bf \ } Let us construct an operator on an odd-dimensional
manifold with a nontrivial fractional part of the $\eta $-invariant. To this
end we apply the cross product of elliptic operators (cf. \cite{AtSi1}).

Let $D_1$ be an elliptic self-adjoint operator on an even-dimensional
manifold $M_1$ with odd symbol:
\[
\sigma \left( D_1\right) \left( x,-\xi \right) =-\sigma \left( D_1\right)
\left( x,\xi \right) ,
\]
and $D_2$ be an elliptic operator on an odd-dimensional $M_2$ with the
symbol satisfying
\begin{equation}
\sigma \left( D_2\right) \left( x,-\xi \right) =\sigma \left( D_2\right)
^{*}\left( x,\xi \right)  \label{alfons}
\end{equation}
(we suppose that $D_2$ is an endomorphism). Denote by $M$ the Cartesian
product of the manifolds $M_1\times M_2.$ Consider the cross product
\[
\left[ \sigma \left( D_1\right) \right] \times \left[ \sigma \left(
D_2\right) \right] \in {\rm {Tor}}K_c^1\left( T^{*}M\right)
\]
of the corresponding elliptic symbols. Here
\[
\left[ \sigma \left( D_1\right) \right] \in K_c^1\left( T^{*}M_1\right)
,\left[ \sigma \left( D_2\right) \right] \in K_c^0\left( T^{*}M_2\right) .
\]
For the (self-adjoint) symbol $\sigma  $ of the product the
following formula is valid
\begin{equation}
\sigma  =\left(
\begin{array}{cc}
\sigma \left( D_1\right) \otimes 1 & 1\otimes \sigma \left( D_2\right) ^{*}
\\
1\otimes \sigma \left( D_2\right) & -\sigma \left( D_1\right) \otimes 1
\end{array}
\right) .  \label{skres}
\end{equation}
Let us compute the linking index of $\left[ \sigma \left( D_1\right) \right]
\times \left[ \sigma \left( D_2\right) \right] $ with the orientation bundle
of $M_1\times M_2.$

\begin{proposition}
\label{hit0}Suppose that $M_2$ is orientable. Then the following equality
holds
\begin{equation}
\left( \left[ \sigma \left( D_1\right) \right] \times \left[ \sigma \left(
D_2\right) \right] \right) \cap \left( 1-\Lambda ^n\left( M^n\right) \right)
=\left[ \sigma \left( D_1\right) \right] \cap \left( 1-\Lambda ^k\left(
M_1\right) \right) {\rm {ind}}D_2,\quad \dim M_1=k.  \label{produ}
\end{equation}
\end{proposition}

\noindent {\em Proof.} Consider the commutative diagram
\[
\begin{array}{ccc}
K_c^0\left( T^{*}M_1,{Q}/{Z}\right) \times K_c^0\left( T^{*}M_2\right) &
\rightarrow & K_c^0\left( T^{*}M,{Q}/{Z}\right) \\[2mm]
\uparrow
\mbox{\footnotesize
$\cap \left( 1-\Lambda ^k\left( M_1\right) \right) \times 1$} &  & \qquad
\uparrow
\mbox{\footnotesize $\cap \left( 1-\Lambda ^k\left( M_1\right)
\right)$} \\[2mm]
{\rm {Tor}}K_c^1\left( T^{*}M_1\right) \times K_c^0\left( T^{*}M_2\right) &
\rightarrow & {\rm {Tor}}K_c^1\left( T^{*}M\right) .
\end{array}
\]
The horizontal mappings here are induced by products in $K$-theory. By
virtue of the orientability of $M_2$ we have
\[
1-\left[ \Lambda ^k\left( M_1\right) \right] =1-\left[ \Lambda ^n\left(
M\right) \right] .
\]
Thence, the desired (\ref{produ}) follows from the last diagram when we
apply the direct image mapping
\[
p_{!}:K_c^0\left( T^{*}M,{Q}/{Z}\right) \rightarrow {Q}/{Z}.
\]
The proposition is proved.

\begin{remark}
{\em Formula (\ref{produ}) is similar to the well-known property of the }$%
\eta ${\em -invariant (see~\cite{APS0}): }$\eta ${\em -invariant of a cross
product is equal to the product of the }$\eta ${\em -invariant of the first
factor and the index of the second operator.}
\end{remark}

Unfortunately, the bundle $L_{+}\left( \sigma  \right) $ is
not an even one. Indeed, the symbol $\sigma  $ satisfies the
equality
\begin{equation}
\alpha ^{*}\sigma  =\left(
\begin{array}{cc}
0 & 1 \\
1 & 0
\end{array}
\right) \sigma  \left(
\begin{array}{cc}
0 & 1 \\
1 & 0
\end{array}
\right) .  \label{involu}
\end{equation}
However, this shows that the spectral subspace transforms according to
\[
\alpha ^{*}L_{+}\left( \sigma  \right) =\left(
\begin{array}{cc}
0 & 1 \\
1 & 0
\end{array}
\right) L_{+}\left( \sigma  \right) .
\]
It turns out that the bundle $L_{+}\left( \sigma  \right) $
is isomorphic to an even bundle.

\begin{proposition}
\label{hit}There exists a vector bundle $L\in {\rm {Vect}}\left(
P^{*}M\right) $ on the projective spaces such that the pull-back $p^{*}L$ to
the cospheres $S^{*}M$ under the projection $p:S^{*}M\rightarrow P^{*}M$ is
isomorphic to $L_{+}\left( \sigma  \right) $.
\end{proposition}

\noindent {\em Proof.} Let us fix a nonsingular vector field on the
odd-dimensional manifold $M$. Consider the corresponding splitting of the
cotangent bundle $T^{*}M=V\oplus 1$ (with respect to some Riemannian
metric). This implies for the cosphere bundle
\[
S^{*}M=S\left( V\oplus 1\right).
\]
The projectivization $P^{*}M$ is diffeomorphic to the ball bundle
\[
BV\subset S\left( V\oplus 1\right)
\]
of $V$ with the identification of antipodal points on the boundary:
\fgr{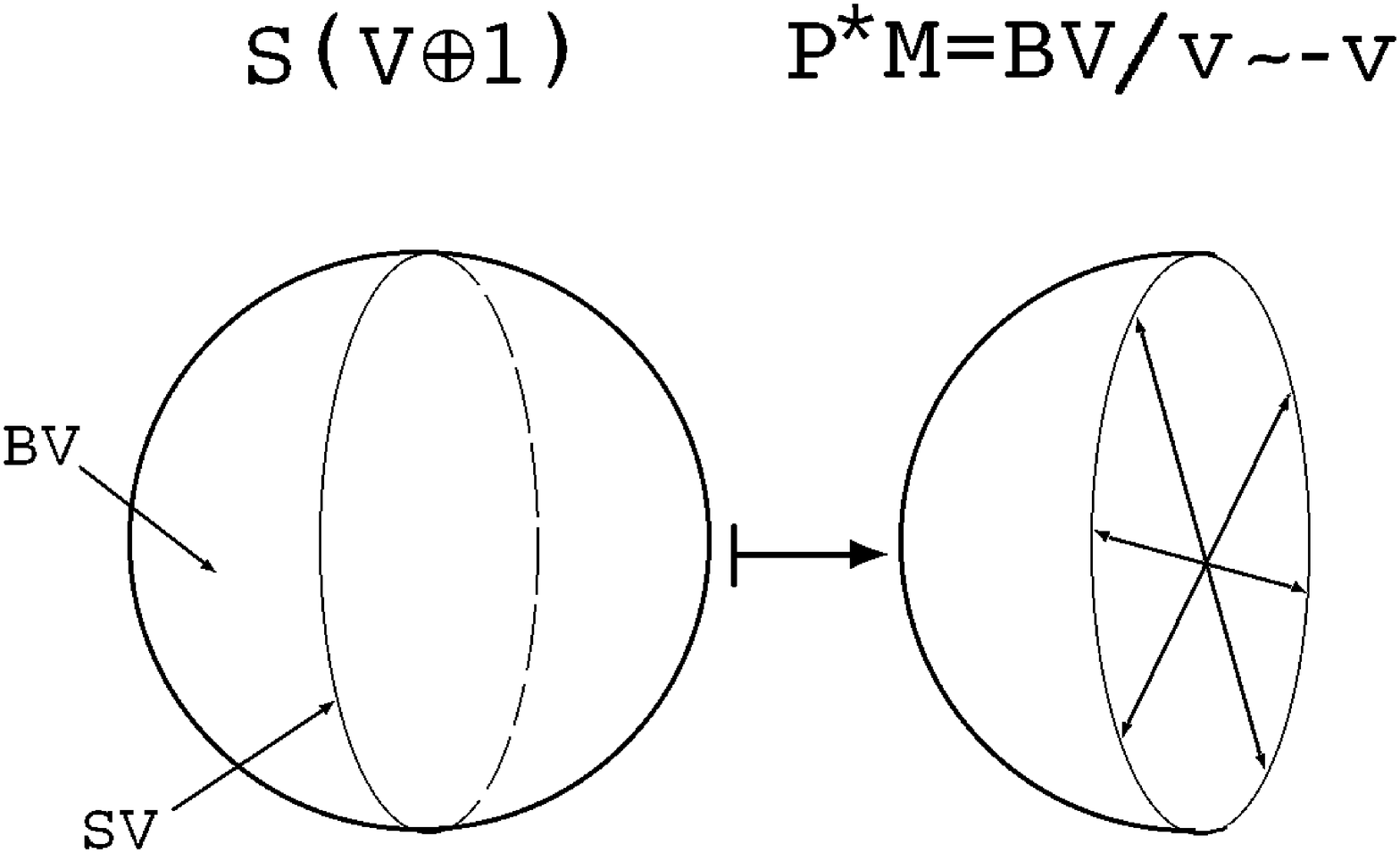}{laboda}{Projective space bundle $P^*M$.}
\[
P^{*}M=BV/\left\{ \left. v\sim-v\right| \,\left| v\right| =1\right\}
\]
(see Fig. 1). Then the bundle $L$ on the projectivization $P^{*}M$ is
constructed similarly from the bundle $L_{+}\left( \sigma
\right) $ on the ball bundle $BV$ by means of the identification of fibers over
antipodes $\pm v$ with respect to the involution $\left(
\begin{array}{cc}
0 & 1 \\
1 & 0
\end{array}
\right) .$ The isomorphism of the pull-back of $L$ to the spheres $S^{*}M$
and the original bundle $L_{+}\left( \sigma  \right) $ is
checked straightforwardly. This proves the proposition.

From the preceding two propositions we immediately obtain a corollary.

\begin{corollary}
The subspace $\widehat{L}$ corresponding to the constructed symbol $L$
satisfies the equality
\begin{equation}
\left\{ 2d\left( \widehat{L}\right) \right\} =\left\{ 2d\left( {\widehat{L}}%
_{+}\left( D_1\right) \right) \right\} {\rm {ind}}D_2.  \label{finform}
\end{equation}
\end{corollary}

{\bf 3.} Let us apply the obtained formula in the following situation. Let $%
M_1={RP}^{2k}$ and $D_1$ be the $pin^c$ Dirac operator from the first
example. As a second factor take the circle $M_2=S^1$ with a
pseudodifferential operator $D_2$ on it with the following principal symbol
\[
\sigma \left( D_2\right) \left( \varphi ,\xi \right) =\left\{
\begin{array}{c}
e^{-i\varphi },\quad \xi=1 , \\
e^{i\varphi },\quad \xi =-1.
\end{array}
\right.
\]
We obtain: ${\rm {ind}}D_2=2$. The dimension functional on the subspace $%
\widehat{L}$ of Proposition \ref{hit} has a nontrivial fractional part, more
precisely
\[
\left\{ 2^{k-1}d\left( \widehat{L}\right) \right\} =\frac 12.
\]
Let us also note one corollary. It gives an answer to the question posed in
\cite{Gil7}.

\begin{corollary}
There exist even-order differential
operators on odd-dimensional manifolds with an arbitrary
dyadic fractional part of the $\eta $-invariant.
\end{corollary}

{\bf 4.} It might seem exotic to consider differential operators of orders
higher that one in index theory. However, there are geometric second-order
operators with interesting spectral properties (see \cite{Gil7}). Let us
also mention that the Hirzebruch operator on an oriented manifold is
equivalent to a second order operator. In particular, the signature of an
orientable manifold is equal to the index of a second order operator. This
observation was used by A.~Connes, D.~Sullivan and N.~Teleman \cite{CST1} to
express the signature of lipschitz and quasiconformal manifolds as an index
of a bounded Fredholm operator.

\vspace{1cm}

\hfill {\em Moscow}

\end{document}